\documentclass[12pt]{article}
\usepackage{amsfonts}
\usepackage{amssymb,amsfonts,amsmath, psfrag,eepic,color}
\usepackage{graphicx, colordvi, cite}

\newtheorem{theo}{Theorem}[section]

\newtheorem{exam}[theo]{Example}
\newtheorem{lem} [theo]{Lemma}

\makeatletter \@addtoreset{equation}{section}
\@addtoreset{theo}{section} \makeatother

\def\pf{\noindent {\textit{Proof.} }}
\def\qed{\hfill \rule{4pt}{7pt}}

\begin{document}
\begin{center}
{\bf \large The Abel-Zeilberger Algorithm}
\vskip 20pt

William Y. C. Chen$^1$, Qing-Hu Hou$^2$, and Hai-Tao Jin$^3$\\[15pt]
Center for Combinatorics, LPMC-TJKLC \\[5pt]
Nankai University, Tianjin 300731, P. R. China\\[5pt]
$^1$chen@nankai.edu.cn, $^2$hou@nankai.edu.cn, $^3$jinht1006@mail.nankai.edu.cn
\vskip 20pt

\it Dedicated to Professor Doron Zeilberger on the occasion of his 60th birthday
\end{center}

\begin{abstract}
We use both Abel's lemma on summation by parts and Zeilberger's
algorithm to find recurrence relations for definite summations. The role
of Abel's lemma can be extended to
the case of linear difference operators
with polynomial coefficients. This approach can be
used to verify and discover identities involving harmonic numbers and derangement numbers.
As examples, we use the Abel-Zeilberger algorithm to  prove the Paule-Schneider identities, the Ap\'{e}ry-Schmidt-Strehl identity, Calkin's identity and some
identities involving Fibonacci numbers.
\end{abstract}



\section{Introduction}

The main idea of this paper is to employ the classical
lemma of Abel and  Zeilberger's algorithm for hypergeometric
sums to verify and to discover various identities
on summations that are not hypergeometric.
Abel's lemma \cite{Abel-1826} on summation
by parts is stated as follows.

\begin{lem}
 For two arbitrary sequences
$\{a_k\}$ and $\{b_k\}$, we have
\[
\sum_{k=m}^{n-1}(a_{k+1}-a_k)b_k=\sum_{k=m}^{n-1}a_{k+1}(b_k-b_{k+1})+a_{n}b_{n}-a_{m}b_{m}.
\]
\end{lem}
For a sequence $\{\tau_k\}$, define the forward
difference operator $\Delta$ by
\[
\Delta\tau_k=\tau_{k+1}-\tau_k.
\]
Then Abel's lemma may be written as
\begin{equation}\label{Abel'lemma}
\sum_{k=m}^{n-1} b_k \Delta a_k =- \sum_{k=m}^{n-1} a_{k+1} \Delta b_k +
a_{n} b_{n} - a_{m}b_m.
\end{equation}

Graham, Knuth and Patashnik \cite{Knuth94} reformulated Abel's lemma in terms of finite calculus to evaluate indefinite sums.
Recently,
Chu \cite{Chu-w-ch2007-1} utilized Able's lemma to prove basic hypergeometric identities
including Bailey's very well poised $_6\psi_6$-series identity by finding pairs $(a_k,b_k)$.  Applying Abel's lemma to the
 pairs $(a_k,b_k)$, one obtains
 contiguous relations for the basic hypergeometric sums. Chen, Chen and Gu
\cite{Chen-2008} presented a systematic approach to finding  pairs $(a_k,b_k)$ by using the $q$-Gosper algorithm.

In this paper, we combine Abel's lemma and Zeilberger's algorithm  to verify and discover identities. Moreover, we use an extended version of Abel's lemma to deal with sums involving holonomic sequences. Let us illustrate our approach by considering identities involving harmonic numbers. The $k$-th harmonic number $H_k$ is given by
\begin{equation}\label{def-H}
H_k = \sum_{j=1}^k \frac{1}{j}.
\end{equation}
Note that by definition, $H_k=0$ whenever $k \le 0$.
Let $f_k$ be a hypergeometric term, i.e., $f_{k+1}/f_k$ is a rational function of $k$. We focus on the summation
\begin{equation}\label{eq-fH}
\sum_{k=m}^{n-1} f_k H_k .
\end{equation}
We can use Gosper's algorithm \cite{Gosper78} to determine whether there exists a hypergeometric term $a_k$ such that $\Delta a_k = f_k$. If such $a_k$ exists,   by Abel's lemma we get
\begin{equation}\label{H-trans}
\sum_{k=m}^{n-1} f_k H_k = \sum_{k=m}^{n-1} H_k \Delta a_k = -\sum_{k=m}^{n-1} \frac{a_{k+1}}{k+1}
+ a_n H_n - a_m H_m.
\end{equation}
Hence we can transform a summation involving harmonic numbers into a hypergeometric summation.
We call such an approach {\it the Abel-Gosper method.}

The same idea applies to the definite summation
\[
S(n) = \sum_{k=0}^n F(n,k) H_k,
\]
where $F(n,k)$ is a proper hypergeometric term such that $F(n,k)=0$ for $k>n$. In this case, we can apply
Zeilberger's algorithm to find a hypergeometric term $G(n,k)$
and polynomials $p_0(n), \dots, p_d(n)$ such that
\[
\sum_{j=0}^d p_j(n)F(n+j,k) =
G(n,k+1)-G(n,k).
\]
Taking
\[
f_k = \sum_{j=0}^d p_j(n)F(n+j,k) \quad \mbox{and} \quad a_k=G(n,k)
\]
in \eqref{H-trans} and summing over $k$ from $0$ to $n+d$, we deduce that
\begin{equation}\label{Abel-Z-H0}
\sum_{j=0}^d p_j(n) S(n+j) = - \sum_{k=0}^{n+d} \frac{G(n,k+1)}{k+1} + G(n,n+d+1) H_{n+d+1}.
\end{equation}
Moreover, we see that the summation on the right hand side of \eqref{Abel-Z-H0} is  again a hypergeometric sum. This approach is called {\it the Abel-Zeilberger method.}

In order to apply the above approach to general holonomic sequences, we extend Abel's lemma by
replacing the operator $\Delta$ with a linear difference operator $L$ of the form
\[
L a_k = r_0(k) a_k + r_1(k) a_{k+1} + \cdots + r_d(k) a_{k+d},
\]
where each $r_i(k)$ is a rational function of $k$.
By combining a variation of Zeilberger's algorithm and the extended Abel's lemma,
we can find a recurrence relation for a sum of the form
\begin{equation}\label{s-ab}
S(n) = \sum_k f(n,k) g(n,k),
\end{equation}
where $f(n,k)$ is hypergeometric and $g(n,k)$ satisfies two recurrence relations
\[
g(n,k+d) = r_1(n,k)g(n,k)+\cdots +r_d(n,k) g(n,k+d-1) + u(n,k),
\]
and
\[
g(n+1,k) = s(n,k) g(n,k) + v(n,k),
\]
where all the coefficients $r_i(n,k)$ and $s(n,k)$ are rational functions and
$u(n,k),v(n,k)$ are hypergeometric terms. The algorithm for finding recurrence relations for the sum \eqref{s-ab} is called the \emph{Abel-Zeilberger algorithm}.

The paper is organized as follows. In Section 2, we give examples to demonstrate that many indefinite sums involving harmonic numbers  can be
reduced to hypergeometric sums by the Abel-Gosper method. Section 3 shows how to
apply the Abel-Zeilberger method to find
recurrence relations of definite sums involving harmonic numbers.  For example, the Paule-Schneider identities fall into this
framework. In Section 4,
 by extending Abel's lemma and using a variation of Zeilberger's algorithm, we present
the Abel-Zeilberger algorithm. The last section provides several examples of the Abel-Zeilberger algorithm including  identities involving Fibonacci numbers and derangement numbers, as well some identities on multiple sums.

\section{The Abel-Gosper method}

In this section, we give several examples to illustrate how to combine Abel's lemma with Gosper's algorithm  to evaluate indefinite summations. We shall focus  on summations involving harmonic numbers.

We begin with a simple example. Consider the sum
\[
S(n) = \sum_{k=1}^{n} H_k.
\]
Expressing $1$ as $\Delta k$, we obtain
\begin{equation} \label{Harmonic-1}
S(n) = \sum_{k=1}^n H_k \Delta k = - \sum_{k=1}^n (k+1) \Delta H_k + (n+1)H_{n+1} - H_1 = (n+1) H_n - n.
\end{equation}

The idea of the above example can be generalized to  indefinite sums of products of polynomials and harmonic numbers. The following result is due to Spie{\ss} \cite{Spieb}. Here we give a derivation based on Abel's lemma.

\begin{theo}\label{Harmonic-3}
Let $H_k$ be the $k$-th harmonic number and $u(n)$ be a polynomial of degree $m$ in $n$. Then
\begin{equation}\label{eq-uH}
\sum_{k=0}^{n} u(k) H_k = p(n)H_n - q(n), \quad n=0,1,2,\ldots,
\end{equation}
where $p(n)$ and $q(n)$ are
both polynomials of degree $m+1$ in $n$. Moreover, $p(n)$ is divisible by $n+1$.
\end{theo}
\pf It is well-known that there exists a polynomial $f(k)$ of degree $m+1$ in $k$
such that
$\Delta f(k) = u(k)$ and the constant term of $f(k)$ is zero. Therefore, we may write $f(k)=k g(k)$, where   $g(k)$ is
 a polynomial of degree $m$. From \eqref{H-trans} it follows that
\[
\sum_{k=0}^{n} u(k) H_k = \sum_{k=0}^{n} \Delta (k g(k)) H_k = -\sum_{k=0}^{n}g(k+1) + (n+1)g(n+1)H_{n+1}.
\]
Since the sum $\sum_{k=0}^n g(k+1)$ is a polynomial of degree $m+1$ in $n$ and
\[
(n+1)g(n+1) H_{n+1} =  (n+1)g(n+1)H_n + g(n+1),
\]
we arrive at \eqref{eq-uH}. \qed

Setting $u(n)=1$ in \eqref{eq-uH}, we obtain \eqref{Harmonic-1}. When $u(n)=n^m$ for $m=1,2,3$, we have
\begin{align}
& \sum_{k=1}^{n}kH_k = \frac{n(n+1)}{2}H_n-\frac{(n-1)n}{4}, \label{Harmonic-p1}\\[6pt]
& \sum_{k=1}^{n}k^2H_k = \frac{n(n+1)(2n+1)}{6}H_n-\frac{(n-1)n(4n+1)}{36},\label{Harmonic-p2}\\[6pt]
& \sum_{k=1}^{n}k^3H_k = \frac{n^2(n+1)^2}{4}H_n-\frac{(n-1)n(n+1)(3n-2)}{48}.\label{Harmonic-p3}
\end{align}

The same idea also applies to the bonus problem 69 proposed by Graham, Knuth and Patashnik \cite[Chapter 6]{Knuth94}.

\begin{exam}
Find a closed form for
\[
\sum_{k=1}^{n} k^2 H_{n+k}.
\]
\end{exam}

The above sum can be rewritten as
\begin{equation}\label{BP69s}
\sum_{k=n+1}^{2n} (k-n)^2 H_k = \sum_{k=1}^{2n} (k-n)^2 H_k - \sum_{k=1}^{n} (k-n)^2 H_k.
\end{equation}
Expanding the summands and applying   formulas \eqref{Harmonic-1}, \eqref{Harmonic-p1}, \eqref{Harmonic-p2}, we obtain that
\begin{equation}\label{BP69a}
\sum_{k=1}^n k^2 H_{n+k} = \frac{n(n+1)(2n+1)}{6}(2H_{2n}-H_n)-\frac{n(n+1)(10n-1)}{36}.
\end{equation}

We may also apply Gosper's algorithm directly to $(k-n)^2$ and get its indefinite sum
\[
\frac{1}{3} k^3 - \frac{2n+1}{2} k^2 + \frac{6n^2+6n+1}{6} k + C.
\]
Setting $C=0$ and using \eqref{H-trans}, we also
obtain \eqref{BP69a}.

We remark that Chyzak \cite{Chyzak2000} and Schneider \cite{Schneider2000, Sch07} have proved \eqref{BP69a} by an extension of Zeilberger's algorithm and Karr's algorithm, respectively.

\begin{exam}
Evaluate the sum
\[
\sum_{k=0}^{n-1} \frac{1}{4^k} {2k \choose k} H_k.
\]
\end{exam}

By Gosper's algorithm, we find that
\[
\Delta \frac{2k}{4^k} {2k \choose k} = \frac{1}{4^k} {2k \choose k}.
\]
Therefore,
\begin{eqnarray*}
\sum_{k=0}^{n-1} \frac{1}{4^k} {2k \choose k} H_k & = &
- 2 \sum_{k=0}^{n-1} \frac{1}{4^{k+1}} {2k+2 \choose k+1} + \frac{2n}{4^n} {2n \choose n} H_n \\[6pt]
&=& 2 - \frac{n+1}{4^n} {2n+2 \choose n+1} + \frac{2n}{4^n} {2n \choose n} H_n.
\end{eqnarray*}

\begin{exam}
We have
\begin{equation}\label{Harmonic-6}
\sum_{k=0}^{n}H_k^2=(n+1)H_n^2-(2n+1)H_n+2n.
\end{equation}
\end{exam}
\pf
Setting $a_k=k$ and $b_k=H_k^2$ in  \eqref{Abel'lemma}, we deduce that
\begin{align*}
\sum_{k=0}^{n} H_k^2 & = -\sum_{k=0}^{n} (k+1) \Delta H_k^2 + (n+1) H_{n+1}^2 \\[6pt]
                     & = -2\sum_{k=0}^{n} H_k - \sum_{k=0}^{n} \frac{1}{k+1} + (n+1)H_{n+1}^2 \\[6pt]
                     & = -2(n+1)H_n+2n -H_{n+1}+(n+1)H_{n+1}^2 \tag*{\mbox{(By \eqref{Harmonic-1})}}\\[6pt]
                     & = (n+1)H_n^2-(2n+1)H_n+2n. \tag*{\qed}
\end{align*}
%\end{exam}

Similarly, by setting $b_k = H_k^3$, we deduce the following identities,
see \cite{Spieb, Chyzak2000},
\begin{align}
\sum_{k=0}^{n}H_k^3&=(n+1)H_n^3-\frac{3}{2}(2n+1)H_n^2+3(2n+1)H_n+\frac{1}{2}H_n^{(2)}-6n,\\
\sum_{k=0}^{n}(2k+1)H_k^3&=(n+1)^2H_n^3-\frac{3}{2}n(n+1)H_n^2+\frac{3n^2+3n+1}{2}H_n-\frac{3}{4}n(n+1).
\end{align}

\section{The Abel-Zeilberger method}

In this section, we use Abel's lemma and Zeilberger's algorithm to find recurrence relations of the definite summation of the form
\[
S(n) = \sum_{k=0}^n F(n,k) H_k,
\]
where $F(n,k)$ is a hypergeometric term in two variables $n$ and $k$. We shall give two examples
to explain the method.
First, we consider an identity due to Chu and De Donno \cite{Chu-w.ch2005}

\begin{exam}
 For $n\geq 0$, we have
\begin{equation}\label{Harmonic-2}
\sum_{k=0}^{n}{n\choose k}^2 H_k=(2H_n-H_{2n}){2n\choose n}.
\end{equation}
\end{exam}

\pf
By applying Zeilberger's algorithm to ${n \choose k}^2$, we obtain a skew
recurrence relation
\[
(n+1) {n+1 \choose k}^2 - 2 (2n+1) {n \choose k}^2 = G(n,k+1) - G(n,k),
\]
where
\[
G(n,k)=(-3-3n+2k){n\choose k-1}^2.
\]
Let
\[
S(n) = \sum_{k=0}^{n}{n\choose k}^2 H_k.
\]
Substituting $F(n,k)={n \choose k}^2$ and $G(n,k)=(-3-3n+2k){n\choose k-1}^2$ into \eqref{Abel-Z-H0}, we find that
\begin{equation}\label{rec-S}
(n+1)S(n+1)-2(2n+1)S(n)= \sum_{k=0}^{n+1} \frac{3n-2k+3}{k+1} {n \choose k}^2 = \frac{4n+1}{n+1}{2n\choose n},
\end{equation}
where the second equality can be justified  by applying Zeilberger's algorithm.

It is easy to verify that \[ R(n)= (2H_n-H_{2n}){2n\choose n}\] satisfies
the same recurrence relation \eqref{rec-S}. Since $S(0)=R(0)=0$, we  get \eqref{Harmonic-2}. This completes the proof. \qed
%\end{exam}

In \cite{Paule-Schneider2003}, Paule and Schneider considered the following summations:
\begin{equation}\label{eq-Tn}
T_n^{(\alpha)} = \sum_{k=0}^{n}(1+ \alpha (n-2k)H_k){n \choose k}^\alpha, \quad \alpha=1,2,\ldots.
\end{equation}
They found closed forms of $T_n^{(\alpha)}$ for $1 \le \alpha \le 4$ and derived recurrence relations of $T_n^{(\alpha)}$ for $5 \le \alpha \le 9$.
As will be seen, we can combine Abel's lemma and Zeilberger's algorithm to deal with the summations $T_{n}^{(\alpha)}$. As an example, let us consider the case $\alpha=3$.

\begin{exam} For $n\geq 0$, we have
\[
T_n^{(3)} = (-1)^n.
\]
\end{exam}
\pf
Let
\[
F(n,k) = (n-2k){n \choose k}^3.
\]
By Zeilberger's algorithm, we find that
\[
F(n,k)+F(n+1,k) = G(n,k+1)-G(n,k),
\]
where
\[
G(n,k) = (2n-k+2){n \choose k-1}^3.
\]
Let
\[
S(n) = T_n^{(3)} = \sum_{k=0}^n (1+3(n-2k)H_k){n \choose k}^3.
\]
By \eqref{Abel-Z-H0}, we deduce that
\[
S(n) + S(n+1) = - 3 \sum_{k=0}^{n+1} \frac{2n-k+1}{k+1} {n \choose k}^3 + \sum_{k=0}^n {n \choose k}^3 + \sum_{k=0}^{n+1} {n+1 \choose k}^3.
\]
By Zeilberger's algorithm, we find that the right hand side, denoted by $R(n)$, satisfies
\[
(n+1)R(n)+(n+2)R(n+1) = 0.
\]
Since $R(0)=0$, we have $R(n)=0$ for $n=0,1,\ldots$.  It is clear that  $S(0)=1$.
 So we get $S(n)=(-1)^n$. This completes the proof.  \qed

Moreover, as a direct consequence of \eqref{Abel-Z-H0}, we have the following property.

\begin{theo}
Let
\[
U_n^{(\alpha)}=\sum_{k=0}^{n}(n-2k){n\choose k}^{\alpha}.
\]
Assume that the minimal recurrence relation for $U_n^{(\alpha)}$ computed by Zeilberger's algorithm is
\[
\sum_{i=0}^{d}p_i(n)U_{n+i}^{(\alpha)}=0.
\]
Then the summation
\begin{equation}\label{SumT}
\sum_{i=0}^{d}p_i(n)T_{n+i}^{(\alpha)}
\end{equation}
is a hypergeometric summation.
\end{theo}

We find that the sum \eqref{SumT} equals zero for $\alpha=1,2,\ldots,9$ and
conjecture that it holds for any nonnegative integer $\alpha$. We note that this conjecture implies the conjecture of  Schneider and Paule  \cite{Paule-Schneider2003}, which says that
 $T_{n+i}^{(\alpha)}$ satisfies the minimal recurrence relation for $U_n^{(\alpha)}$ computed by Zeilberger's algorithm.

\section{The Abel-Zeilberger algorithm}

In this section, we give a description of the Abel-Zeilberger algorithm.
Notice that in the applications of Abel's lemma given in previous sections, the main idea lies in the fact that $\Delta H_k$ is a hypergeometric term. In fact, there are other  sequences satisfying similar properties that lead us to
consider an extension of Abel's lemma.

Let $\{a_k\}$ be an arbitrary sequence. We consider a linear operator $L$ of the form
\[
L a_k = r_0(k) a_k + r_1(k) a_{k+1} + \cdots + r_d(k) a_{k+d},
\]
where each $r_j(k)$ is a rational function of $k$. We associate the operator $L$ with a dual operator $L^*$ defined by
\[
L^* a_k  = r_0(k)a_k + r_1(k-1)a_{k-1} + \cdots + r_d(k-d) a_{k-d}.
\]
 In the above notation, Abel's lemma can be extended to the following form.

\begin{lem} For two arbitrary sequences
$\{a_k\}$ and $\{b_k\}$, we have
\begin{equation}\label{Generalized Abel's lemma}
\sum_{k=m}^{n-1} L^* a_k \cdot b_k = \sum_{k=m}^{n-1} a_k \cdot Lb_k
- T(n) + T(m),
\end{equation}
where
\begin{equation}\label{def-T}
T(k) = \sum_{i=1}^d \sum_{j=1}^i r_i(k-j)a_{k-j}b_{k+i-j}.
\end{equation}
\end{lem}

\pf It is easy to verify that
\begin{align*}
\sum_{k=m}^{n-1} L^* a_k \cdot b_k &=
\sum_{k=m}^{n-1} \sum_{i=0}^{d} r_i(k-i)a_{k-i} b_{k}\\[6pt]
&= \sum_{i=0}^{d}\sum_{k=m-i}^{n-1-i}r_i(k) a_k b_{k+i}\\[6pt]
&= \sum_{k=m}^{n-1}r_0(k)a_{k}b_{k}+\sum_{i=1}^{d}\left[
\sum_{k=m}^{n-1}+\sum_{k=m-i}^{m-1}-\sum_{k=n-i}^{n-1}\right]r_{i}(k)a_k b_{k+i}\\[6pt]
&= \sum_{k=m}^{n-1}a_{k}\cdot
L b_{k}+\sum_{i=1}^{d}\sum_{k=m-i}^{m-1}r_{i}(k) a_k b_{k+i} - \sum_{i=1}^{d} \sum_{k=n-i}^{n-1} r_{i}(k) a_k b_{k+i}.  \tag*{\qed}
\end{align*}

Let $f(n,k)$ be a bivariate hypergeometric term and $g(n,k)$ be a bivariate function. We aim to find a linear recurrence relation for the definite sum
\[
S(n) = \sum_{k=m}^{\ell}  f(n,k) g(n,k).
\]
Suppose that there exist rational functions $r_j(n,k)$ such that
\[
L g(n,k) = \sum_{j=0}^d r_j(n,k) g(n,k+j)
\]
is a bivariate hypergeometric term. We shall try to find polynomials $p_i(n)$, which are independent of $k$ and not all zero, together with hypergeometric terms $a(n,k)$ and $w(n,k)$ such that
\begin{equation}\label{paw}
\sum_{i=0}^I p_i(n) f(n+i, k) g(n+i, k) = g(n,k) L^* a(n,k) + w(n,k).
\end{equation}
Summing \eqref{paw} over $k$ and applying the extended Abel's lemma, we deduce that
\begin{align}
\sum_{i=0}^I p_i(n) S(n+i) & =  \sum_{i=0}^I \sum_{k=m}^{\ell}  p_i(n) f(n+i, k) g(n+i, k) \nonumber \\
&= \sum_{k=m}^{\ell} ( g(n,k) L^* a(n,k) + w(n,k) ) \nonumber \\
&= \sum_{k=m}^{\ell} a(n,k) L g(n,k) + \sum_{k=m}^\ell w(n,k) - T(\ell+1) + T(m), \label{nonhomrec}
\end{align}
where $T(k)$ is given by \eqref{def-T}. Notice that in the last expression, the two summands are hypergeometric. By Zeilberger's algorithm, the two sums satisfy linear recurrence relations, which lead to a non-homogenous linear recurrence relation for $S(n)$.

To solve equation \eqref{paw},
we impose the condition  that $g(n,k)$ satisfies
the relation
\[
g(n+1,k) = s(n,k)g(n,k) + v(n,k),
\]
where $s(n,k)$ is a rational function and $v(n,k)$ is a hypergeometric term. By induction, it is
easy to show that
there exist rational functions $s_i(n,k)$ and hypergeometric terms $v_i(n,k)$ such that
\begin{equation}\label{sivi}
g(n+i,k) = s_i(n,k) g(n,k) + v_i(n,k).
\end{equation}

Now we can solve the following equation for $p_i(n)$ and $a(n,k)$ by a variation of Zeilberger's algorithm
\begin{equation}\label{pa}
\sum_{i=0}^I p_i(n) f(n+i, k) s_i(n,k) = L^* a(n,k).
\end{equation}
It can be seen that $a(n,k)$ is similar to $f(n,k)$, that is,
\[
R(n,k)= \frac{a(n,k)}{f(n,k)}
\]
is a rational function of $n$ and $k$. Hence \eqref{pa} is equivalent to
\begin{equation} \label{eq-R}
\sum_{i=0}^I p_i(n) \frac{f(n+i, k)}{f(n,k)} s_i(n,k) = \sum_{j=0}^d \frac{f(n,k-j)}{f(n,k)} r_j(n,k-j) R(n,k-j).
\end{equation}
Since $f(n,k)$ is hypergeometric, both $f(n+i,k)/f(n,k)$ and $f(n,k-j)/f(n,k)$ are rational functions. Therefore, \eqref{eq-R} is a non-homogenous linear recurrence equation on $R(n,k)$ with parameters $p_i(n)$, which can be solved by Abramov's algorithm \cite{Abramov}.

Once we find a solution $(p_0(n), \ldots, p_I(n), a(n,k))$ to \eqref{pa}, it is easy to check that
\[
(p_0(n), \ldots, p_I(n), a(n,k), w(n,k))
\]
is a solution to \eqref{paw}, where
\begin{equation}\label{eq-w}
w(n,k) = \sum_{i=0}^I p_i(n) f(n+i, k) v_i(n,k).
\end{equation}

In summary,  the Abel-Zeilberger algorithm can be described as follows.

\begin{description}
\item
{\it Input:} a hypergeometric term $f(n,k)$ and a term $g(n,k)$ satisfying two recurrence relations
\begin{equation}\label{rec-1}
g(n,k+d) = r_1(n,k)g(n,k)+\cdots +r_d(n,k) g(n,k+d-1) + u(n,k),
\end{equation}
and
\begin{equation}\label{rec-2}
g(n+1,k) = s(n,k) g(n,k) + v(n,k),
\end{equation}
where $r_i(n,k)$ and $s(n,k)$ are rational functions, and
$u(n,k)$ and $v(n,k)$ are hypergeometric terms.

\item
{\it Output:} polynomials $p_i(n)$ that are independent of $k$, two hypergeometric terms $t_1(n,k)$, $t_2(n,k)$ and a term $T(k)$ satisfying
\[
\sum_{i=0}^I p_i(n) S(n+i) = \sum_{k=m}^{\ell} t_1(n,k) + \sum_{k=m}^{\ell} t_2(n,k) - T(\ell+1) + T(m),
\]
where
\[
S(n) = \sum_{k=m}^\ell f(n,k)g(n,k).
\]
\end{description}

The algorithm consists of the following steps.

\noindent Initially, we set $I=0$.

\noindent {Step 1.}
For $0 \le i \le I$, compute the rational functions $s_i(n,k)$ and the hypergeometric terms $v_i(n,k)$ defined by \eqref{sivi} by using the recurrence relations
\begin{align*}
& s_{i+1}(n,k)=s(n+i,k) s_i(n,k), \\[6pt] & v_{i+1}(n,k) = s(n+i,k)v_i(n,k) + v(n+i,k),
\end{align*}
with the initial values $s_0(n,k)=1$ and $v_0(n,k)=0$.

\noindent {Step 2.}
Let
\begin{equation}\label{def-L}
L g(n,k) = -r_1(n,k) g(n,k) - \cdots -r_d(n,k) g(n,k+d-1) + g(n,k+d).
\end{equation}
According to \eqref{eq-R}, construct an equation on $p_i(n)$ and $R(n,k)=a(n,k)/f(n,k)$. That is, compute polynomials
\[
P_j(n,k), \, 0 \le j \le d \quad \mbox{and} \quad Q_i(n,k), \, 0 \le i \le I
\]
such that
\begin{equation} \label{PRQ}
\sum_{j=0}^d P_j(n,k) R(n,k-j) = \sum_{i=0}^I p_i(n) Q_i(n,k).
\end{equation}

\noindent {Step 3.}
Solve equation \eqref{PRQ} for $R(n,k)$ and $p_i(n), 0 \le i \le I$ by using Abramov's algorithm.
If all the polynomials $p_i(n)$ are zeros, then we increase $I$ by one and repeat steps 1--3.

\noindent {Step 4.}
Compute $w(n,k)$ based on \eqref{eq-w}. For $L$ given by \eqref{def-L}, compute $T(k)$ according to \eqref{def-T}.
Finally, set
\[
t_1(n,k) = a(n,k) u(n,k) \quad \mbox{and} \quad  t_2(n,k) = w(n,k).
\]
Then $(p_i(n), t_1(n,k), t_2(n,k), T(k))$ is the desired output.

\section{Examples}
In this section, we provide several examples to
compute summations by using the Abel-Zeilberger algorithm.

We first consider the case when $g(n,k)$ is independent of $n$ in  the Abel-Zeilberger algorithm as described in the previous section, i.e., $s(n,k)=1$ and $v(n,k)=0$ in \eqref{rec-2}. In this case, we have $s_i(n,k)=1$ and $v_i(n,k)=0$ so that $t_2(n,k)=0$.

Let $F_n$ be the $n$-th \emph{Fibonacci number} which is defined by the recurrence relation
\[
F_{n+2}=F_{n+1}+F_n,
\]
 for $ n \ge 0$, with initial values $F_0=0, F_1 = 1$. Employing the Abel-Zeilberger algorithm, we can prove the following identities on the product of binomial coefficients and Fibonacci numbers, see \cite{Benjamin, Vajda1989}.
\begin{exam}\label{Fibonacci-1} We have
\begin{align}
& \sum_{k=0}^n {n\choose k}F_k = F_{2n}, \label{F1}\\[6pt]
& \sum_{k=0}^n (-1)^k {n\choose k} F_k = -F_{n}, \\[6pt]
& \sum_{k=0}^n {n\choose k}F_{3k}=2^nF_{2n}, \\
& \sum_{k=0}^n {n\choose k}F_{4k}=3^nF_{2n}.
\end{align}
\end{exam}

\pf For equation \eqref{F1}, taking $f(n,k)={n \choose k}$ and
$g(n,k)=F_k$ as the input of the Abel-Zeilberger algorithm, we
obtain
\[
p_0(n)=1,\ p_1(n)=-3,\ p_2(n)=1, \quad t_1(n,k)=t_2(n,k)=0,
\]
and
\[
T(k) = \big( (2k-n-3)F_k + (n+2-k)F_{k+1} \big) \frac{n!}{(n+2-k)!(k-1)!}.
\]
We see that $T(n+3)=T(0)=0$. Therefore, the summation
\[
S(n) = \sum_{k=0}^n {n \choose k} F_k
\]
satisfies the recurrence relation
\[
S(n)-3S(n+1)+S(n+2)=0.
\]
Since $F_{2n}$ satisfies the same recurrence relation with initial values $F_0 = 0$ and $F_2 = 1$,
we deduce that $S(n)=F_{2n}$.

The other three identities can be proved in the same fashion.
 The detailed arguments are omitted. This completes the proof.
\qed
%\end{exam}

It is not difficult to see that the Abel-Zeilberger algorithm could be used to verify identities involving general $C$-finite sequences. For instance, suppose that $\{G_k\}_{k \ge 0}$ satisfies a recurrence relation
\[
G_{k+2}=b G_{k+1} + c G_k, \quad k\geq 0,
\]
where $b$ and $c$ are constants. Then the Abel-Zeilberger algorithm generates the recurrence relation
\[
(b+1-c)S(n)-(b+2)S(n+1)+S(n+2)=0
\]
for the summation
\[
S(n)=\sum_{k=0}^{n} {n\choose k} G_k.
\]

The next example involves the $n$-th \emph{derangement number}
$D_n$ as given by
\[
D_n = n! \sum_{k=0}^{n}\frac{(-1)^k}{k!}.
\]
Using the method of MacMahon's partition analysis,
 Andrews and Paule \cite{G.E.Andrews1999}
 the following identity
 on $D_n$.
We shall  give a derivation
 by applying the Abel-Zeilberger
algorithm.

\begin{exam}

\begin{equation}\label{G.E.Andrews-1999}
\sum_{j\geq0} {k \choose j} \frac{(k+n-j)!}{(k+N-j)!}
D_{k+N-j} = (-1)^n \sum_{j\geq0} (-1)^j {n\choose
j} (j+k)!,
\end{equation}
for  $n\geq N\geq n-k$.
\end{exam}

\pf
Substituting $k+N-j$ for $j$,  the left hand side of \eqref{G.E.Andrews-1999} can be rewritten as
\[
S(N,n) = \sum_{j=N}^{k+N} {k\choose j-N}\frac{(n+j-N)!}{j!} D_j.
\]
Because of the recurrence relation
\[
D_n = n D_{n-1} + (-1)^n
\]
for $D_n$, we may put
\[
f(N,j)={k\choose j-N}\frac{(n+j-N)!}{j!}
 \quad \mbox{and} \quad g(N,j)=D_j
\]
as the input of the Abel-Zeilberger algorithm.
Then we obtain
\[
S(N,n)-S(N+1,n) = \sum_{j\geq N}(-1)^j{k\choose
j-N}\frac{(n+j-N)!}{(j+1)!} .
\]
Denote the right hand side by $G(k)$.
 By Zeilberger's algorithm,
we find that for $ k \ge 0$,
\[
(1+N-n+k)G(k)-(2+k+N)G(k+1)=0,
\]
which implies that
$G(k) = 0$ for  $k \ge n-N \ge 0$. Thus we get
$S(N,n)=S(N+1,n)$. In particular,
\[
S(N,n)=S(n,n)=\sum_{j\geq n}{k\choose j-n} D_j.
\]
Applying the Abel-Zeilberger algorithm to
\[
f(n,j) = {k \choose j-n} \quad \mbox{and} \quad g(n,j)=D_j,
\]
we find
\[
(n+1)S(n,n)+(n+k+1)S(n+1,n+1)-S(n+2,n+2) = \sum_{j\geq n}(-1)^j {k+1\choose
j-n} = 0.
\]
By Zeilberger's algorithm, we see that
 the right hand side of \eqref{G.E.Andrews-1999}
 satisfies the same recurrence relation.
 Finally, from the identity
\[
\sum_{j=0}^k {k \choose j} D_j = k!,
\]
we deduce that $S(0,0)=k!$ and $S(1,1)=(k+1)!-k!$,
which coincides with the initial values of
 the right hand side of \eqref{G.E.Andrews-1999}.
  Thus \eqref{G.E.Andrews-1999}
   holds for $n \ge N \ge n-k$. \qed
%\end{exam}

The following identity was found by Schmidt and has been proved in several
ways, see \cite{Strehl-94, Wegschaider1997}.

\begin{exam}
[The Ap\'{e}ry-Schmidt-Strehl Identity]
 For $n\geq 0$, we have
\begin{equation}\label{Wei-P71-Identity}
\sum_{k}\sum_{j}{n\choose k}{n+k\choose k}{k\choose
j}^3=\sum_{k}{n\choose k}^2{n+k\choose k}^2.
\end{equation}
\end{exam}

\pf
Let
\[
f(n,k) = {n\choose k} {n+k\choose k},  \quad
g(n,k) =\sum_{j=0}^k {k\choose j}^3,
\]
and
\[
S(n) = \sum_{k=0}^n f(n,k)g(n,k).
\]
By Zeilberger's algorithm, $g(n,k)$ is annihilated by the operator
\[
-8(k+1)^2-(7k^2+21k+16)K+(k+2)^2K^2.
\]
where $K$ denotes the shift operator on $k$. Now applying the Abel-Zeilberger algorithm to
$f(n,k)$ and $g(n,k)$, we obtain
\[
S(n)-\frac{(3+2n)(39+51n+17n^2)}{(n+1)^3}S(n+1)+\frac{(n+2)^3}{(n+1)^3}S(n+2)=0,
\]
Meanwhile, using Zeilberger's algorithm, we find that the Ap\'{e}ry
numbers \[ A(n) = \sum_{k=0}^n {n\choose k}^2{n+k\choose k}^2\]  satisfy
the same recurrence relation.
Finally, by comparing the initial values, we arrive at \eqref{Wei-P71-Identity}. \qed
%\end{exam}

To conclude this paper, we consider the following summations,
\[
S_n^{(\alpha)}=\sum_{k=0}^{n}\left(\sum_{j=0}^{k}{n\choose j}\right)^{\alpha}.
\]
For $\alpha=1,2$ and $3$, closed forms for $S_n^{(\alpha)}$ have been derived by Andrews and Paule \cite{G.E.Andrews1999} by using the
method of MacMahon's partition analysis.
It is easy to see that these formulas can be
derived by using the Abel-Zeilberger algorithm.

\begin{exam}\label{ex-Calkin}
We have
\begin{equation}\label{Calkin's identity-1}
S_n^{(1)} = \sum_{k=0}^{n}\sum_{j=0}^{k}{n\choose j} = n2^{n-1}+2^n,
\end{equation}
and
\begin{equation}\label{Calkin's identity-2}
S_n^{(2)} = \sum_{k=0}^{n}\left(\sum_{j=0}^{k}{n\choose
j}\right)^2=\left(\frac{n}{2}+1\right)2^{2n}-\frac{n}{2}{2n\choose n}.
\end{equation}
\end{exam}
\pf Let
\[
f(n,k)=1, \quad  g(n,k)=\sum_{j=0}^{k}{n\choose j}.
\]
It is clear that
\begin{equation}\label{fg-rec}
g(n,k+1)=g(n,k)+{n\choose k+1} \quad \mbox{and} \quad
g(n+1,k)=2g(n,k)-{n\choose k}.
\end{equation}
Applying the Abel-Zeilberger algorithm to $f(n,k)$ and $g(n,k)$, we find that
\[
S_n^{(1)} = \sum_{k=0}^n (C-k){n \choose k+1} - T(n+1) + T(0),
\]
where $T(k) = (C+1-k) g(n,k)$ and $C$ is a constant. Setting $C=-1$, we get $T(0)=0$. Thus
\begin{align*}
S_n^{(1)}&=-\sum_{k=0}^{n}(k+1){n\choose k+1}+(n+1)g(n,n+1)\\[6pt]
         &=-n\sum_{k=0}^{n}{n-1\choose k}+(n+1)\sum_{j=0}^{n+1}{n\choose j}\\[6pt]
         &=-n2^{n-1}+(n+1)2^n\\[6pt]
         &=n2^{n-1}+2^n.
\end{align*}

Now we consider the evaluation of $S_n^{(2)}$. Let $f(n,k), g(n,k)$ be given as above, and let $h(n,k)=g(n,k)^2$. From \eqref{fg-rec}, we see that
\begin{align*}
& h(n,k+1)=h(n,k)+u(n,k), \\[6pt]
& h(n+1,k) = 4h(n,k) + v(n,k),
\end{align*}
where
\begin{align*}
& u(n,k) = 2{n\choose k+1}g(n,k)+{n\choose k+1}^2,
\\[6pt]
& v(n,k) = -4{n\choose k}g(n,k)+{n\choose k}^2.
\end{align*}
It should be mentioned that there is actually
no need to impose the condition for $u(n,k)$ and $v(n,k)$ to be hypergeometric in the Abel-Zeilberger algorithm. Therefore, we can still apply the Abel-Zeilberger algorithm to $f(n,k)$ and $h(n,k)$ to deduce that
\[
S_n^{(2)} =\sum_{k=0}^{n}(C-k) \left( 2 {n\choose
k+1}g(n,k)+{n\choose
k+1}^2 \right) + (n-C) h(n,n+1)+(C+1)h(n,0),
\]
where $C$ is a  constant. Setting $C=n/2-1$ and applying the Abel-Zeilberger algorithm, we find that
\[
\sum_{k=0}^{n}(n/2-1-k) {n\choose
k+1}g(n,k) = \sum_{k=0}^{n} \frac{-n+k+1}{2} {n\choose k+1}^2 - \frac{n}{2} g(n,0).
\]
Hence
\begin{align*}
S_n^{(2)} &=\sum_{k=0}^{n}(-n+k+1){n\choose k+1}^2-n +
       \sum_{k=0}^{n} \left( \frac{n}{2}-1-k \right) {n\choose k+1}^2 +
       \left( \frac{n}{2}+1 \right) 2^{2n} + \frac{n}{2}\\[6pt]
    &=-\frac{n}{2}\sum_{k=0}^{n}{n\choose k+1}^2 + \left( \frac{n}{2}+1 \right) 2^{2n}-\frac{n}{2}\\[6pt]
    &=-\frac{n}{2}\left({2n\choose n}-1\right)+ \left( \frac{n}{2}+1 \right) 2^{2n}-\frac{n}{2}\\[6pt]
    &=\left( \frac{n}{2}+1 \right)2^{2n} - \frac{n}{2}{2n\choose n}. \tag*{\qed}
\end{align*}

Using the same argument, we can deduce Calkin's identity \cite{Calkin}
\begin{equation}
S_n^{(3)} = n2^{3n-1}+2^{3n}-3n2^{n-2}{2n\choose n}.\label{Calkin's
identity}
\end{equation}

\noindent\textbf{Acknowledgments.}  We wish to thank Professor Lu Yang for valuable comments.
 This work was supported by the
973 Project, the PCSIRT Project of the Ministry of Education, and the National Science
Foundation of China.

%---------------------------------------------------------------

\end{document}